\documentclass{amsart}
\usepackage{amsmath, amsfonts, amssymb, latexsym}
\newtheorem{theorem}{Theorem}
\newtheorem{lemma}{Lemma}

\newtheorem{assumption}[lemma]{Assumption}

\newtheorem{remark}[lemma]{Remark}

\begin{document}
\title{Nonexistence of solutions in $(0,1)$ for K-P-P-type equations
 for all $d\ge 1$}
\author{J\'{a}nos Engl\"{a}nder and P\'{e}ter L. Simon}
\address{Department of Statistics and Applied Probability\\
University of California, Santa Barbara\\
CA 93106-3110, USA\\ and ELTE, Department of Applied Analysis\\
H-1117 Budapest, P\'{a}zm\'{a}ny P\'{e}ter S\'{e}t\'{a}ny 1/C,
Hungary.} \email{englander@pstat.ucsb.edu, simonp@cs.elte.hu}
\urladdr{http://www.pstat.ucsb.edu/faculty/englander,
http://www.cs.elte.hu/~simonp} \keywords{KPP-equation, semilinear
elliptic equations, positive bounded solutions, branching
Brownian-motion} \subjclass{Primary: 35J60, 35J65; Secondary: 60J80}
\date{\today}

\begin{abstract}
Consider the KPP-type equation of the form $\Delta u+f(u)=0$, where
$f:[0,1] \to \mathbb R_{+}$ is a concave function. We prove for
arbitrary dimensions that there is no solution bounded in $(0,1)$.
The significance of this result from the point of view of
probability theory is also discussed.
\end{abstract}
\maketitle

\section{Introduction and main result}
In this article we will investigate certain semilinear elliptic
equations of the form $\Delta u+f(u)=0$. Our assumption on the
nonlinear term $f(u)$ is as follows.
\begin{assumption}\rm
We assume that $f:[0,1] \to \mathbb R$ is
\begin{enumerate}
\item[(i)]  continuous , \item[(ii)]  positive in $(0,1)$ and
\item [(iii)] $z\mapsto f(z)/z$ is strictly
decreasing.$\hfill\diamond$
\end{enumerate}
\end{assumption}
\noindent Consider now the Kolmogorov Petrovskii Piscunov-type
(KPP) equation

\begin{eqnarray}
\Delta u+f(u)=0  \label{eqn} \\
0<u<1 , \  \mathrm{in}\ \mathbb R^d.\label{cond}
\end{eqnarray}
\begin{theorem}\label{mainthm}
Problem (\ref{eqn})-(\ref{cond}) has no solution for $d\ge 1$.
\end{theorem}
\medskip
Semilinear elliptic equations of the form (\ref{eqn}) have been
widely studied. We mention here only two reviews \cite{Ouyang,
Tang}, where the exact number of positive solutions with different
nonlinearities are studied. In \cite{Ouyang} the differential
equation is considered on a bounded domain, in \cite{Tang} the
equation is studied in the whole space $\mathbb R$, however, it is
subject to the boundary condition $u \to 0$ as $|x| \to \infty$. The
case of concave $f$ has also been studied by several authors. In
\cite{Brezis} the assumption on $f$ is similar to ours, however, the
problem is given on a bounded domain with Dirichlet boundary
condition. In that paper the existence and uniqueness of the
positive solution is proved. Castro et al. studied the case of
concave nonlinearities in a series of papers, see e.g.
\cite{Castro1, Castro2}. In these works the problem is given on a
bounded domain with Dirichlet boundary condition. A generalized
logistic equation, with $f(u)=mu-qu^p$ is studied in
\cite{Hernandez} on a bounded domain with Dirichlet boundary
condition again.

Summarizing, we can say that our equation (\ref{eqn}) has been
widely studied, however, in the papers where it is considered in the
whole space $\mathbb R$, it is always subject to the boundary
condition $u \to 0$ as $|x| \to \infty$. In these publications the
aim is to determine the exact number of the so-called fast and slow
decay solutions. Hence according to the authors knowledge there is
no result available concerning problem (\ref{eqn})-(\ref{cond})
under the assumptions given on $f$.

\begin{remark}[Low dimensions]\rm Our theorem can be proved very
 easily for $d\le 2$.
To see this, recall that $\Delta$ is a so-called {\it critical}
operator in $\mathbb R^d$ when $d=1,2$. Second order elliptic
operators $L$ with no zeroth order term are classified as being
{\it subcritical} or critical according to whether the operator
possesses or does not possess a minimal positive Green's function.
In probabilistic terms criticality/subcriticality is captured by
the {\it recurrence/transience} of the corresponding diffusion
process (see Chapter 4 in \cite{P}).

Another equivalent condition for $L$  to be critical is that all
positive functions $h$ that are superharmonic (i.e. $Lh\le 0$) are
in fact harmonic (i.e. $Lh\equiv 0$). (See again Chapter 4 in
\cite{P})

Now, observe that (\ref{eqn})-(\ref{cond}) and the positivity of
$f$ on $(0,1)$ implies
\begin{eqnarray}
\Delta u=-f(u)< 0  \label{eqn.lowdim}\  \mathrm{in}\ \mathbb
R^d.\label{cond.lowdim}
\end{eqnarray}
By the above criterion for critical operators, this  is impossible
in dimension one or two.$\hfill\diamond$

\end{remark}

The most important model case is the classical KPP equation, when
\begin{equation}\label{specialNL}
 f(u):=\beta u(1-u)
\end{equation}
with $\beta>0$. (In fact this particular nonlinearity is
intimately related to the distribution of a {\it branching
Brownian motion}; see more on the subject in the next paragraph.)
Here we present a proof of this result which is valid basically
for concave functions. In fact, (iii) of Assumption 1  is related
to the concaveness of the function.

The connection between the KPP equation and branching Brownian
motion has already been discovered by McKean  --- it first appeared
in the classic work \cite {McK1975,McK1976}.

Let $Z=(Z(t))_{t\geq 0}$ be the $d$-dimensional binary branching
Brownian motion  with a spatially and temporally constant branching
rate $\beta>0$. The informal description of this process is as
follows. A single particle starts at the origin, performs a Brownian
motion on $\mathbb R^d$, after a mean--$1/\beta$ exponential time
dies and produces two offspring, the two offspring perform
independent Brownian motions from their birth location, die and
produce two offspring after independent mean--$1/\beta$ exponential
times, etc. Think of $Z(t)$ as the subset of $\mathbb R^d$
indicating the locations of the particles $z_1^t,...,z^{N_t}_t$
alive at time $t$ (where $N_t$ denote the number of particles at
$t$). Write $P_{x}$ to denote the law of $Z$ when the initial
particle starts at $x$. The natural filtration is denoted by
$\{\mathcal {F}_t,\ t\ge 0\}$.

Then, as is well known (see e.g. Chapter 1 in \cite{D02}), the law
of the process can be described via its Laplace functional as
follows. If $f$ is a positive measurable function, then
\begin{equation}\label{Laplace.func}
E_x\exp \left(-\sum_{i=1}^{N_t} f(z_i^t)\right)=1-u(x,t),
\end{equation}
where $u$ solves the initial value problem
\begin{eqnarray}
&&\dot{u}=\frac{1}{2}\Delta u+f(u) \  \mathrm{in}\
\mathbb R^d\times\mathbb R_+  \label{par.eqn} \\
&&u(\cdot,0)=1-e^{-f(\cdot)} \  \mathrm{in}\
\mathbb R^d\nonumber\label{IC}\\
&&0\le u\le 1 \  \mathrm{in}\ \mathbb R^d\times\mathbb
R_+,\nonumber\label{same.cond}
\end{eqnarray} with $f$ from
(\ref{specialNL}).

Choosing appropriate (sequences of) $f$'s one can then express the
probabilities of various events from $A\in \mathcal{F}_t$, for
$t>0$, in terms of the function $u$ in (\ref{par.eqn}). Letting
$t\to\infty$ then routinely leads to results stating that if $A\in
\mathcal{F}_{\infty}$ denotes a certain event then the function
$u(x):=P_x(A)$ is either constant ($=0$ or $=1$), or it must solve
(\ref{eqn})-(\ref{cond}). So, if one knows our main theorem then
it immediately follows that all those events are {\it trivial}
(that is, their probability is either zero or one).

Equations of the type (\ref{eqn})-(\ref{cond}) frequently appear
when one studies certain `natural' martingales associated with
branching Brownian motion (see e.g. \cite{EK}).

Note that if $\beta>0$ is replaced by a smooth nonnegative function
$\beta(\cdot)$ that does not vanish everywhere, then this
corresponds to having {\it spatially dependent} branching rate for
the branching Brownian motion.

It would be desirable therefore to investigate whether our main
theorem can be generalized for such $\beta$'s.
\section{Proof of the theorem}

The proof is based on two ideas: the application of the semilinear
elliptic maximum principle, which is generalized here fore concave
functions, and a comparison between the semilinear and the linear
problems. Using these two ideas we will show that the {\it minimal
positive solution} of (\ref{eqn}) is $u_{\min}\equiv 1$, hence
(\ref{eqn}) has no solution satisfying (\ref{cond}).

First we state and prove a semilinear maximum principle.
The results in this form is a
generalization of \cite[Proposition 7.1]{EP99}
for the particular case when the elliptic
operator is $L=\Delta$.

\begin{lemma}[Semilinear elliptic maximum principle]\label{emp}
Let $f:[0,\infty ) \to \mathbb R$ be a continuous function, for
which Assumption 1(iii)  holds. Let  $D\subset \Bbb{R}^{d}$ be a
bounded domain with smooth boundary. If $v_{i}\in C^2(D)\cap
C(\bar D)$\ satisfy $v_{i}>0$ in $D$, $\Delta v_{i}+f( v_{i})=0,\
in\ D$ for $i=1,2$, and $v_{1}\ge v_{2}$ on $\partial D$, then
$v_{1}\ge v_{2}$ in $\bar D$.
\end{lemma}




\medskip
\noindent\underline{Proof}: The function $w:=v_1-v_2$ satisfies
\begin{equation}
\Delta w + f(v_1)-f(v_2)=0 . \label{3}
\end{equation}
We show that $w\ge 0$  in $D$. Suppose to the contrary that there
exists a point $y\in D$ where $w$ is negative. Let $\Omega_0:=\{x\in
D\mid w(x)<0\}$. Let $\Omega$ be the connected component of
$\Omega_0$ containing $y$. Since $w \geq 0$ on $\partial D$, one has
$\Omega \subset \subset D$ and
\begin{equation}
w<0 \ \mbox{ in } \Omega  \qquad  w=0 \ \mbox{ in } \partial \Omega . \label{4}
\end{equation}
 Let us multiply the
equation $\Delta v_{1}+f( v_{1})=0$ by $w$ and equation (\ref{3})
by $v_1$, then subtract the second equation from the first, and
integrate on $\Omega$. Using that $w=v_1-v_2$ one obtains
\begin{equation}
I+II:=\int\limits_{\Omega} (w \Delta v_1 - v_1 \Delta w) +
\int\limits_{\Omega} (v_1 f(v_2) -v_2 f(v_1)) =0. \label{5}
\end{equation}
Using Green's second identity and that $w=0$ in $\partial \Omega$
along with the fact that $\partial_{\nu} w \geq 0$ on $\partial
\Omega$, we obtain
$$
I= -\int\limits_{\partial \Omega} v_1
\partial_{\nu} w \leq 0,
$$
where $\nu$ denotes the unit outward normal to $\partial \Omega$.
Furthermore, since  $v_1<v_2$ in $\Omega$, using (iii) of
Assumption 1, we have
 that also $II<0$:
$$
v_1 f(v_2) - v_2 f(v_1)= v_1 v_2 \left[ \frac{f(v_2)}{v_2} -
\frac{f(v_1)}{v_1} \right] < 0.
$$
It follows that the left hand side of (\ref{5}) is negative, while
its right hand side is zero. This contradiction proves that in fact
$w\ge 0$ in $D.\hfill \blacksquare$

\begin{remark}[Spatially dependent $f$'s]\rm
One can similarly prove the analogous more general result for the
case, when $f:D\times [0,\infty ) \to \mathbb R$ is continuous in
$u$ and bounded in $x$, and $u \mapsto f(x,u)/u$ is strictly
decreasing.$\hfill\diamond$
\end{remark}
Let $f:[0,1] \to \mathbb R$ be a continuous function which is
positive in $(0,1)$. Based on ideas in \cite{KS} and using the
comparison between the linear and the semilinear equations,  we
prove the following lemma.

\begin{lemma}[Radially symmetric solutions] \label{lemma2}
Assume in addition that $f$ satisfies $\liminf\limits_{z
\downarrow 0} \frac{f(z)}{z}>0$ (this is automatically satisfied
under Assumption 1(iii)). Then for any $y\in \mathbb R^d$ and
$p\in (0,1)$ there exists a ball $\Omega:= B_R(y)$ (with some
$R>0$) and a radially symmetric $C^2$ function $v:\Omega \to
\mathbb R$ such that
\begin{eqnarray}
\Delta v + f(v)&=& 0 \nonumber \\
v>0 \ \mbox{ in } \Omega, & & \ v=0 \ \mbox{ in } \partial\Omega, \ v(y)=p . \nonumber
\end{eqnarray}
\end{lemma}

\medskip
\noindent\underline{Proof}: We show the existence of  a radially
symmetric solution of the form $v(x)=V(|x-y|)$. Let $V\in
C^2([0,\infty))$ be the solution of the initial value problem
\begin{eqnarray}
(r^{d-1}V'(r))'+r^{d-1}f(V(r))=0 \label{6} \\
V(0)=p, \ V'(0)=0 .
\end{eqnarray}
Writing $\Delta$ in polar coordinates, one sees that it is
sufficient to prove that there exists an $R>0$ such that $V(R)=0$
and $V(r)>0$ for all $r\in [0,R)$.

To this end, consider the {\it linear} initial value problem
\begin{eqnarray}
(r^{d-1}W'(r))'+r^{d-1}m W(r)=0 \label{7} \\
W(0)=p, \ W'(0)=0 ,
\end{eqnarray}
where $m>0$ is chosen so that $f(u)>mu$ holds for all $u\in (0,p)$.
(Our assumptions on $f$ guarantee the existence of such an $m$.) It
is known that $W$ has a first root, which we denote by $\rho$. Note
that in this case $-m$ is the first eigenvalue of the Laplacian on
the ball $B_{\rho}$. We now show that $V$ has a root in $(0,\rho]$.
In order to do so let us multiply (\ref{7}) by $V$ and (\ref{6}) by
$W$, then subtract one equation from the other, and finally,
integrate on $[0,\rho]$. We obtain
\begin{eqnarray}
&&I+II:=\int\limits_0^{\rho} [(r^{d-1}W'(r))'V(r) -
(r^{d-1}V'(r))'W(r)]\, \mbox{d}r\nonumber\\ &&\ \ \ \ \ \  \ \ \
+\int\limits_0^{\rho} r^{d-1}[mW(r)V(r) - W(r)f(V(r))]\, \mbox{d}r =
0 . \label{8}
\end{eqnarray}
Suppose now that $V$ has no root in $(0,\rho]$. Then, integrating by
parts, $ I = \rho^{d-1}W'(\rho)V(\rho)<0. $

Next, observe that by integrating (\ref{6}), one gets $V'(r)<0$
(i.e. $V$ is decreasing). Hence $V(r)<p$, yielding
$mV(r)-f(V(r))<0.$ Therefore $II$, and thus the whole left hand side
of (\ref{8}) are negative; contradiction. This contradiction proves
that $V$ in fact has a root in $(0,\rho].\hfill\blacksquare$

\begin{remark}[Spatially dependent f's]\rm
When $f$ depends also on $x$, our method breaks down as it is no
longer possible to use ordinary differential equations to show the
existence of a solution attaining a value close to one at a given
point.

There is one easy case though: it is immediately seen that if
there exists a $g(u)$, with $f(x,u)\ge g(u)$ and $g(u)$ satisfies
the conditions of Theorem \ref{mainthm}, then Theorem
\ref{mainthm} remains valid for $f(x,u)$ as well.

Indeed, we know that $u_{min}\ge 1$, where $u_{min}$ is the minimal
positive solution for the semilinear equation with $g$. Recall (see
e.g. \cite{EP99, EP03}) that one way of constructing the minimal
positive solution is as follows. One takes  large balls $B_R(0)$,
and positive solutions  with zero boundary condition on these balls
(in our case  we know from \cite{KS} that there exist such positive
solutions for arbitrarily large $R$'s), and finally, lets
$R\to\infty$; using the monotonicity in $R$ that follows from the
semilinear elliptic maximum principle (Lemma \ref{emp}), the
limiting function exists and positive. It is standard to prove that
it solves the equation on the whole space, and  by Lemma \ref{emp}
again it must be the {\it minimal} such solution.

Now suppose that $0<v$ solves the semilinear equation with $f(x,u)$.
Then  $v$ is a {\it supersolution}: $0\ge \Delta v+ g(v)$; hence by
the above construction of $u_{min}$ and by an obvious modification
of the proof of Lemma \ref{emp}, $v\ge u_{min}\ge 1$.

The general case is harder. For example, when
$f(x,u):=\beta(x)(u-u^2)$ and $\beta$ is a smooth nonnegative
bounded function, the mere existence of positive solutions on large
balls is no problem  as long as the generalized principal eigenvalue
of $\Delta+\beta$ on $\mathbb R^d$ is positive. (The method in
\cite{Pconst}, pp. 26-27 goes through for $f(x,u):=\beta(x)(u-u^2)$
even though $\beta$ is constant in \cite{Pconst}.) The problematic
part is to show that the solution is large at the center of the
ball. $\hfill \diamond$
\end{remark}

\medskip \noindent\underline{Proof of Theorem \ref{mainthm}}:
Suppose that problem (\ref{eqn})-(\ref{cond}) has a solution. Choose
an arbitrary point $y\in \mathbb R^d$ and an arbitrary number $p\in
(0,1)$. Note that by Assumption 1, $f$ satisfies the conditions of
Lemma \ref{lemma2} and consider  the ball $B_R(y)$ and the radially
symmetric function $v$ on it, which are guaranteed by Lemma
\ref{lemma2}. We can apply Lemma \ref{emp} with $D=B_R(y)$, $v_1=u$
and $v_2=v$ and obtain that $u\geq v$. In particular then, $u(y)\geq
v(y)=p$. Since $y$ and $p$ were arbitrary, we obtain that $u\geq 1$,
in contradiction with (\ref{cond}). Consequently,
(\ref{eqn})-(\ref{cond}) has no solution.$\hfill\blacksquare$

\end{document}